\documentclass[1p, authoryear]{elsarticle}
\usepackage{amsmath}
\usepackage{graphicx}
\usepackage{epstopdf}
\usepackage{natbib}
\usepackage{esint}
\usepackage{dirtytalk}
\usepackage[normalem]{ulem}
\usepackage[colorlinks=true, linkcolor=blue, anchorcolor=red, citecolor=blue, filecolor=red, urlcolor=blue, pdfauthor=author]{hyperref}
\usepackage[font=scriptsize,caption=false]{subfig}
\usepackage{alphalph}

\usepackage{multirow}
\usepackage[dvipsnames]{xcolor}
\usepackage{mwe}
\usepackage{makecell}
\usepackage{hhline}
\usepackage{float}
\usepackage{bookmark}
\usepackage{tabu}
\usepackage{relsize}
\usepackage{upgreek}
\usepackage{bm}
\usepackage{uri} 

\journal{Engineering Computations}

\begin{document}

\begin{frontmatter}
	
	\title{\textbf{\normalsize The official version of this paper can be downloaded at \href{https://www.emerald.com/insight/content/doi/10.1108/EC-01-2021-0052/full/html}{https://doi.org/10.1108/EC-01-2021-0052}.} \\ \vspace{1em}
	Revisiting Non-Convexity in Topology Optimization of Compliance Minimization Problems}
	\author{Mohamed Abdelhamid}
	\author{Aleksander Czekanski\corref{mycorrespondingauthor}}
	\address{Department of Mechanical Engineering, Lassonde School of Engineering, York University, Toronto, Ontario, Canada}
	\cortext[mycorrespondingauthor]{Corresponding author}
	\ead{alex.czekanski@lassonde.yorku.ca}
\begin{abstract} 
	
\noindent \textbf{Purpose --} This is an attempt to better bridge the gap between the mathematical and the engineering/physical aspects of the topic. We trace the different sources of non-convexification in the context of topology optimization problems starting from domain discretization, passing through penalization for discreteness and effects of filtering methods, and end with a note on continuation methods.

\noindent \textbf{Design/Methodology/Approach --} Starting from the global optimum of the compliance minimization problem, we employ analytical tools to investigate how intermediate density penalization affects the convexity of the problem, the potential penalization-like effects of various filtering techniques, how continuation methods can be used to approach the global optimum, and how the initial guess has some weight in determining the final optimum.

\noindent \textbf{Findings --} The non-convexification effects of the penalization of intermediate density elements simply overshadows any other type of non-convexification introduced into the problem, mainly due to its severity and locality. Continuation methods are strongly recommended to overcome the problem of local minima, albeit its step and convergence criteria are left to the user depending on the type of application.

\noindent \textbf{Originality/Value --} In this article, we present a comprehensive treatment of the sources of non-convexity in density-based topology optimization problems, with a focus on linear elastic compliance minimization. We put special emphasis on the potential penalization-like effects of various filtering techniques through a detailed mathematical treatment.


\noindent \textbf{Paper Type} Research paper

\end{abstract}
\end{frontmatter}

\section{Introduction}
Since the introduction of the \textbf{s}olid \textbf{i}sotropic \textbf{m}aterial with \textbf{p}enalization (SIMP) method in the seminal paper by \citet{Bendsoe1989}, material interpolation methods have become one of the most active research areas in engineering optimization. Although the origin of almost all density-based approaches lies in linear elasticity, they have been successfully extended to even more complicated single and multiphysics fields. In a general sense, the optimal topology of a problem refers to the location and number of holes such that an objective function is extremized.

To simplify the numerical implementation of topology optimization formulations, the normalized density is usually taken to be element-wise constant rather than the less popular node-wise approach. This way the design variables can be taken outside the integral of the elemental stiffness matrices (cf. \citet[p.~68]{Sigmund1998} and \citet[p.~1419]{Petersson1998}). A discrete formulation, where the normalized density is only allowed to take discrete \( 0/1 \) values, proves to be extremely difficult to solve numerically \citep{Stolpe2007}. Hence relaxation of the design space is introduced, namely continuous\footnote{\, To avoid confusion, \textit{continuous} is used to describe the opposite of discrete design variable, while \textit{continuum} is used to describe the opposite of discretized geometrical spatial domains.} variation of the design variables, where intermediate density elements between void and material states are allowed to exist during the optimization process \citep{Bendsoe1999}. In order to obtain practical manufacturable designs, intermediate density elements have to be eliminated from the final design. Penalization, in different forms, is introduced to increase the cost (i.e., reduce the benefit) of intermediate elements to force the solution to converge towards discrete \(0/1 \) designs. The concept of continuous design variables with penalization constitutes the essence of material interpolation methods.

It is a widely known fact that any approach that enforces discrete \( 0/1 \) solutions is inherently non-convex. A rather unfavorable consequence of this non-convexity is the obscureness of the global optimum. In other words, the converged solution of non-global optimization approaches would be one of the huge pool of local minima the problem possesses. This is mainly due to the fact that in a non-convex function, the global optimum can't be identified using certain mathematical conditions, as opposed to the local optimum that can be identified through the local behavior of the function, namely its gradient and Hessian. Hence, in a non-convex problem, it is impossible to calculate, or even claim, the global optimum without an exhaustive search of the entire design space \citep[p.~138-139]{Arora1995}. To complicate the problem further, global optimization approaches are extremely computationally extensive, and, so far, prove unable to handle the typically massive number of design variables in practical topology optimization problems \citep[p.~74]{Sigmund1998}.

A major numerical instability that was early recognized as a direct consequence of the non-convexification of the objective function is the \textit{local minima} problem. It mainly refers to the problem of obtaining different solutions for the same \textbf{f}inite \textbf{e}lement (FE) discretization upon choosing different algorithmic parameters \citep[p.~74]{Sigmund1998}.
The converged solution of an optimization problem is mainly determined by three parameters\footnote{\, Strictly speaking, this statement only applies to optimization problems with a non-changing objective function, so it doesn't strictly apply when continuation methods are used (i.e., increasing the penalization factor). The use and effect of continuation methods is discussed in depth in a later section.}; the initial guess, the optimization direction, and the optimization \say{speed}. The latter two parameters are characteristics of the optimization algorithm/solver in use \citep{Bendsoe1999}. The optimization direction is determined by the gradients\footnote{\, The discussion in this study is solely dedicated to gradient-based optimization, since the usefulness of non-gradient based methods is still unclear as was discussed in the forum article by  \citet{Sigmund2011}.} of the objective function and constraints with respect to the design variables (i.e., sensitivity analysis), while the optimization speed is determined by how aggressive the algorithm is (e.g., how much change in the design variables is allowed per iteration).
So far, numerous studies attempted to tackle the last two parameters (i.e., optimization direction and speed) through focusing the attention on the type of the optimization algorithm and its controls. As for the initial guess, \citet[p.~74]{Sigmund1998} mentioned in their discussion on local minima that even with the use of artificial schemes to convexify the problem and lead to reproducible designs, small variations in initial parameters such as geometry of design domains could lead to substantial changes in the converged design. They also mentioned the necessity for using continuation methods, and cited some literature on different continuation procedures.
\citet[p.~298]{Cardoso2003} mentioned in their proposal of a new general mesh independent filter, how the initial guess could have an effect on the converged solution, and proposed the use of continuation methods to mitigate this effect. \citet{aremu2010suitability} observed how the converged solution of a random initial guess is often different, both numerically and topologically, from the typical result of an even distribution initial guess. \citet{Yan2018} emphasized the significance of using carefully selected initial guesses to circumvent the existence of many local minima in volume-to-point heat conduction optimization problems, and suggested utilizing very slow continuation approaches or additional design constraints to overcome the problem.

Even though the final desired solution in a topology optimization problem is ideally pure discrete with no intermediate density elements, the discrete form of the problem is almost impossible to challenge. Hence, in a discussion on global optimality, the logical path to follow is to start from the relaxed convex form of the compliance minimization problem that has an easily-identifiable global optimum, that is the variable thickness sheet problem. Starting from this global optimum, we investigate how penalization affects the convexity of the problem, how continuation methods can be used to approach the global optimum, and how the initial guess has some weight in determining the final optimum.

The rest of the paper is organized as follows: \textbf{Section \ref{sec_probform}} outlines the general formulation for the compliance minimization problem of interest, \textbf{Section \ref{sec_var_thick_sh_prob}} discusses the relaxed form known as the variable thickness sheet problem and its properties and significance, \textbf{Section \ref{sec_penal}} details the concept of penalization, its sources and effects on the convexity of the problem, and how it can be remedied through continuation methods, \textbf{Section \ref{sec_challenges}} discusses some computational methods relevant to non-convexity and the problem of local minima, and finally the main conclusions of the article are discussed in \textbf{Section \ref{sec_conc}}.

\section{Problem Formulation}
\label{sec_probform}
In this section, for the sake of clarity, we detail the general problem formulation utilized in all the upcoming discussions. It is readily observable that most topology optimization problems cannot be solved in the continuum formulation and have to be discretized into a number of finite elements. Hence, the problem formulation in this study is stated directly in the discretized form of the design domain. A more detailed discussion on the effect of domain discretization is included in Section \ref{sec_var_thick_sh_prob}.

So far, there has been a number of different density-based approaches (also termed \textit{material interpolation methods}) that adopt penalization to approximate discrete \( 0/1 \) solutions. It's currently well-known that different interpolation methods behave differently based on the characteristics of the problem of interest. In other words, the solution takes a different trajectory for each interpolation method \citep{Stolpe2001a}. However, we believe their non-convexification effects from the global optimality point of view share enough similarities to warrant a generalization for this discussion. Hence, the discussions in the current study will focus on the most popular approach, the SIMP formulation, with references to other methods where appropriate. The typical compliance minimization problem using the \textit{modified} SIMP approach can be stated as follows:
\begin{equation}
\begin{aligned}
\underset{\bm{\uprho}}{\text{min:}} & \qquad c(\bm{\uprho}) = \mathbf{U}^T \mathbf{K} \mathbf{U} . \\
\text{subject to:}                     & \qquad \mathbf{E}(\bm{\uprho}) = E_\text{min}+ \bm{\uprho}^p (E_0-E_\text{min}) , \\[4pt]
                                       & \qquad \mathbf{K} \mathbf{U} = \mathbf{F} , \\[4pt]
                                       & \qquad \mathbf{0} \leq \bm{\uprho} \leq \bm{1} , \\[4pt]
                                       & \qquad \sum_{i=1}^{N} V_i/V_0 \leq V_f .
\label{eqn_prob_form}
\end{aligned}
\end{equation}
where \( c(\boldsymbol{\uprho}) \) is the compliance of the structure, \( \boldsymbol{\uprho} \) is the vector of normalized elemental densities, \( \boldsymbol{\mathrm{U}} \) is the global displacement vector, \( \boldsymbol{\mathrm{K}} \) is the global stiffness matrix, \( \boldsymbol{\mathrm{E}}(\boldsymbol{\uprho}) \) is the vector of elemental elastic moduli, \(E_\text{min}\) is the elastic modulus of the void material, \( p \) is the penalization factor, \( E_0 \) is the elastic modulus of the solid material, \( \boldsymbol{\mathrm{F}} \) is the global force vector, \( V_i \) is the volume of element \( i \), \( N \) is the total number of finite elements, \( V_0 \) is the original design domain volume, and \( V_f \) is the prescribed volume fraction.

It's worth noting that the above formulation implements the \textit{nested approach}. The nested formulation, as opposed to the \textit{\textbf{s}imultaneous \textbf{a}nalysis \textbf{a}nd \textbf{d}esign} (SAND) formulation, is often more suitable for iterative solvers, plus it enables the use of off-the-shelf commercial FE solvers for the equilibrium equations. On the other hand, the simultaneous formulation introduces an additional source of non-convexity that stems from the governing equations. In addition, since SAND treats both the design and the state variables as unknowns, large scale problems have been historically disadvantaged using this approach (cf. \citet[p.~24 \& 222]{Bendsoe2004} and \citet[p.~85]{Christensen2009}). However, recent work on SAND have revived the interest in this approach's ability to solve large scale problems \citep{Munro2017, Munro2017a, Munro2018}. 

Furthermore, the discussions in this study are limited to 2D problems without any loss of generality. The conclusions made should easily extend to 3D problems since there isn't any inherent difference in the behavior of 2D vs. 3D problems from the optimization point of view.

\section{The Variable Thickness Sheet Problem}
\label{sec_var_thick_sh_prob}
The topology optimization problem in the continuum state is ill-posed, as the global optimum would be the configuration in the limit of an infinite number of infinitesimal microscopic holes. In other words, there is always a better global optimum by decreasing the size and increasing the number of these microscopic holes, which means the solution set is not closed. From a mathematical sense, the correct term to use would be a global \textit{infimum}, not a minimum, meaning it can be approached but never attained.

In order to fully comprehend the effect of domain discretization on convexity, it's worthwhile to revisit the definition of \textit{convex functions}. \citet[p.~98]{Bazaraa2006} define a function to be convex if a line segment drawn between any two points on its graph falls entirely on or above the graph. This definition has two implications; one on the \textit{domain} of the function and another on its \textit{codomain}. Considering first the \textit{domain} (i.e., the set of solution points), it has to form a convex set. Meaning it has to be continuous and its boundaries must not curve into the set. As for the \textit{codomain}, for a function \(f(\mathbf{x})\) with \(S\) being a nonempty set in \(R^N\), convexity requires that the function's epigraph (i.e., the set of points on or above its graph) defined as:
\begin{equation}
\{(\mathbf{x},y): \, \mathbf{x} \in S, \, y \in R, \, y \geq f(\mathbf{x})\}.
\label{eq_epi}
\end{equation}

\noindent be a convex set itself. The first implication relates to domain discretization while the second relates to penalization.

\begin{figure*}
    \centering
    \subfloat[]{\label{f1_a} \includegraphics[width=0.49\textwidth]{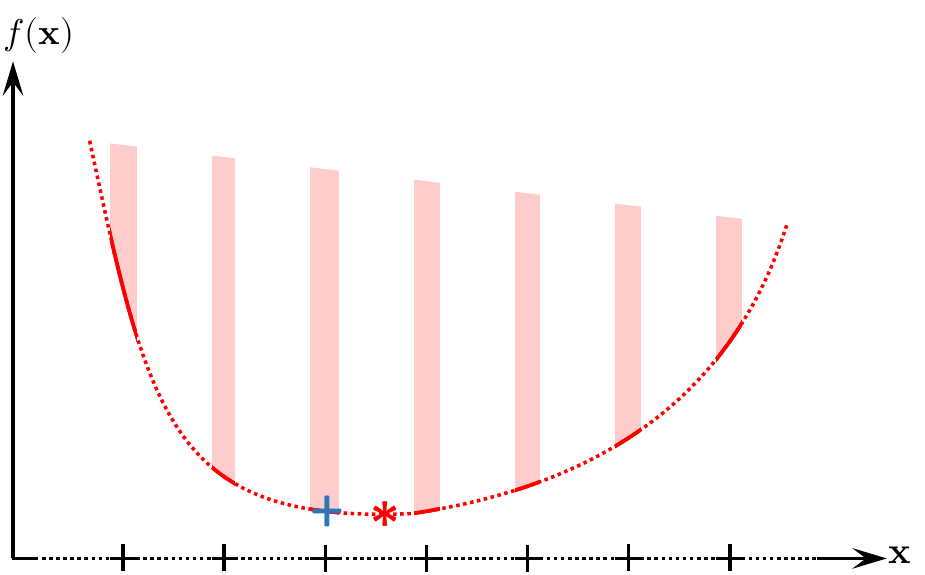}}
    \subfloat[]{\label{f1_b} \includegraphics[width=0.49\textwidth]{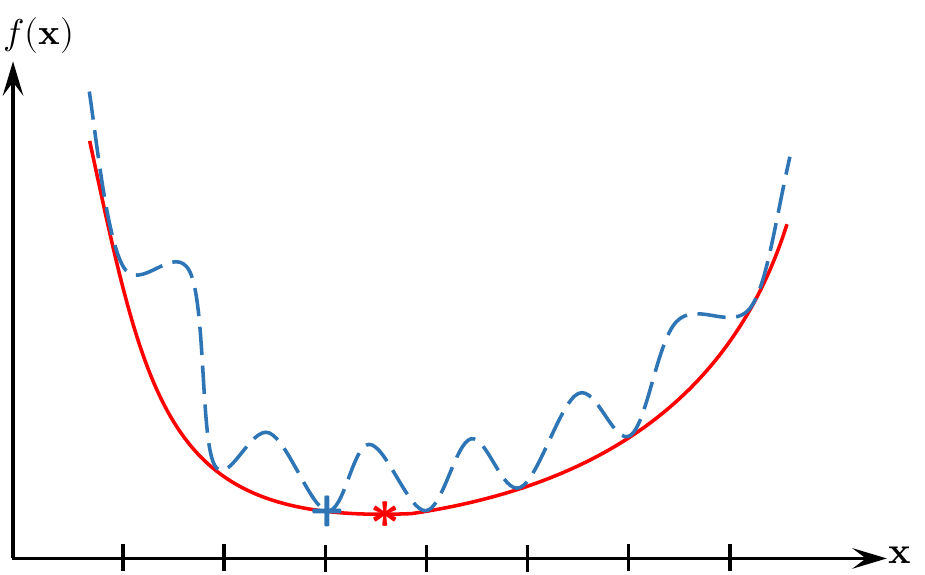}}
    \caption{Non-convexification of a convex function by: \textbf{a. Domain Discretization -} the dotted lines represent the solution points removed due to domain discretization, while the light red areas represent the function's epigraph after discretization. \textbf{b. Penalization -} the original convex function is shown in continuous red, while the penalized non-convex version is shown in dashed blue. In both cases, the global minimum is marked before non-convexification by a red asterisk, and after non-convexification by a blue plus symbol. The tick marks on the horizontal axis represent the desirable design points, whereas the undesirable points are represented by the values in between.}
    \label{f_non-convexification}
\end{figure*}

Returning to our discussion on the variable thickness sheet problem, for obvious practical reasons, the size and the number of the holes have to be within certain limits, typically a lower limit for the size and an upper limit for the number of the holes. Once these limits are imposed, non-convexification is introduced. This is precisely the effect of transforming the continuum domain into a discretized one; the number and size of the potential holes is limited as a result of this discretization \citep{Kohn1986a,Kohn1986b,Kohn1986c}. From a mathematical perspective, domain discretization\footnote{\, Although the main domain discretization technique typically employed in topology optimization problems is the finite element method, other techniques have been successfully adopted such as the finite cell method in \cite{Parvizian2012}. The discussion in this section should apply to any domain discretization technique without loss of generality.} causes non-convexity by \say{poking holes} in the original convex set of solution points, that is by removing design points with intermediate density elements from the set of solutions points. This concept is illustrated in Fig. \ref{f1_a}. Nonetheless, domain discretization remedies the lack of closedness of the solution set, simply by limiting the available options of distributing a certain amount of material within a discretized domain. This is in addition to the obvious need for domain discretization since all practical problems are impossible to solve in the continuum state.

To alleviate the non-convexity triggered by domain discretization, \say{relaxation} has to be introduced. Relaxation in this context refers to expanding the solution space to include the potential solution/design points eliminated by domain discretization (in one form or another). In topology optimization, two main relaxation techniques exist; \textbf{(i)} homogenization using \textit{optimal}\footnote{\, \textit{Optimal} in the sense that the microstructure is capable of attaining the optimal energy bounds; that is utilizing a certain amount of material to the fullest.} microstructures such as rank-2 laminates, and \textbf{(ii)} non-penalized material interpolation methods which is the variable thickness sheet problem. Even though relaxation creates a well-posed problem, it is a double-edged approach. For instance, consider the ideal case of relaxation introduced by using optimal microstructures such as rank-2 laminates. On the one hand, relaxation restores the convexity and produces the closest solution possible to the original continuum problem (for specific displacement and density fields). On the other hand, it also means that the topological features of the optimal design are inherently exhibited on the micro scale, not on the macro scale. The ingenuity of the homogenization-based approach introduced by \citet{Bendsoe1988} is that they introduced \textit{sub-optimal} microstructures, square elements with square holes, that inherently include penalization to force the optimal topology to appear on the macro scale  \citep[p.~210]{Allaire1993a}. One way to think about domain discretization combined with relaxation is as follows; \textbf{(i)} the original problem in the continuum formulation consists of \textit{a solution space of infinite dimensions} \(R^\infty\) with \textit{a discrete design variable} \(\{0,1\}\) for each dimension, \textbf{(ii)} domain discretization with relaxation transforms the problem into \textit{a solution space of finite dimensions} \(R^N\) with \textit{a continuous design variable} \([0,1]\) for each dimension. The beauty of the transformed formulation is that it's numerically bounded, albeit still non-practical.

Given its simplicity of numerical implementation compared to the homogenized rank-2 laminates, from this point forward, our focus in discussion would be on the variable thickness sheet problem without loss of generality. The variable thickness sheet problem is, in its general form, a flat sheet loaded within its plane and the thickness is allowed to vary at every point within the sheet. The characteristic feature of the variable thickness sheet problem is its linear dependence of stiffness and volume on the sheet thickness. This leads easily to existence of solutions for the compliance problem, even without any limitations on the solution space. The variable thickness sheet problem has a physical meaning in its 3D form, however, it cannot be considered a pure 2D problem. Instead, the 2D form of this problem is achieved through utilizing the Voigt upper bound on stiffness \citep{Swan1997}, which still retains the well-behaviour with respect to existence of solutions for the compliance objective function \citep{Bendsoe1999}. Work on this problem was done as early as the paper by \citet{Prager1968}. The first finite element numerical solution to this problem was done by \citet{Rossow1973}. In \citep{Petersson1999}, the author proved the convergence and uniqueness of solution for the finite element discretization of the variable thickness sheet problem. Notably the author proposed two assumptions to frame his argument: \textbf{(i)} \textit{regularity} of the stress field, which refers to the smoothness of the stress field over the problem domain. In engineering terms, this assumption translates to using conforming finite elements, in which case compatibility of the mesh is achieved and monotonic finite element convergence is attainable (cf. \citet[p.~377]{Bathe2006finite}), and \textbf{(ii)} \textit{biaxiality} of the load; this is mainly to avoid the non-uniqueness of solutions of uniaxial loading problems where a thick strip of material would be equivalent in load-carrying efficiency to several thin strips of material of the same volume\footnote{\, This point is elaborated in detail in \citet[p.~70]{Sigmund1998} in their discussion on non-uniqueness of solutions.}. Both assumptions make perfect sense from the engineering point of view. A noteworthy remark about the work by \citet{Petersson1999} is that he had to use higher order elements to avoid the numerical instability of checkerboarding. This issue is discussed in detail in the next subsection.
\par \vspace{1.5em}

\textbf{A Side Note on Checkerboarding} \par \vspace{1em}
Ever since the earlier attempts in topology optimization, the numerical instability of checkerboarding was observed. This issue was first observed in the work of \citet{Cheng1981} during their investigation of thin elastic Kirchhoff plates. Moreover, this issue is similar to previous observations in shape optimization where \textit{zigzag} boundary conditions, though not optimal, were obtained when the boundary nodal coordinates were used directly as design variables \citep{Zhou2001a, Haftka1986}. The current consensus among the research community is that checkerboarding is related to how the displacement and density fields are discretized, even though the continuum problem might be well-behaved \citep{Brezzi1991}. It occurs, in principle, through the intricate interaction of two factors: \textbf{(i)} the kinematic constraints imposed on deformation depending on the type and formulation of the finite elements in use. For instance, bilinear four-node quad elements have the inherent constraint that element edges are always linear. Such constraints produce numerically-induced artificial high stiffness, that would otherwise vanish due to stress singularities at the checkerboard corner connections \citep[p.~42]{Diaz1995}, and \textbf{(ii)} penalization of intermediate density elements either through material interpolation methods or through homogenization with sub-optimal microstructures. The second factor simply means that utilizing the same amount of material is more efficient if used in solid \say{black/white} rather than in penalized intermediate \say{gray} elements.

Checkerboarding was studied in detail almost simultaneously in the pioneering work by \citet{Diaz1995} and \citet{Jog1996}. The first work took more of an engineering approach and limited their study to element-wise constant density formulation, which is the most popular implementation. The authors computed the effective properties of checkerboard arrangements for a volume fraction of \(0.5\) while imposing kinematic constraints on the displacements to ensure consistency with typical finite element approximations. These checkerboard patches were considered for two types of finite elements; bilinear four-node and biquadratic nine-node quad elements. Upon comparing the strain energy densities of these checkerboard patches with a variety of optimal and sub-optimal microstructures, they concluded that in the case of the four-node bilinear quad elements, checkerboarding is prone to appear for any microstructure representation; that is using homogenization or material interpolation methods. Moreover, in a shocking contradiction to typical intuition, they proved mathematically that four-node quad elements can also exhibit checkerboarding within the relaxed state; that is using non-penalized material interpolation methods or using homogenization with optimal microstructures. As for the nine-node biquadratic quad elements, the authors concluded that checkerboarding can occur, albeit only for high penalization factors depending on the value of Poisson's ratio.

On the other hand, \citet{Jog1996} took more of a mathematical approach where they studied the stability of checkerboards through considering the topology optimization problem as a mixed variational problem similar to those used to model incompressible elasticity and Stokes flow. They considered a variety of interpolations for both the displacement and the density fields. Their analysis included developing an incremental form of the variational problem that could be used within an iterative solution procedure for the nonlinear problem. In the relaxed state, the stability of the optimization problem was related to a condition, in the proximity of the optimal solution, on the rank of a matrix whose negative columns represent the pseudoload vectors corresponding to the elements of the incremental density variables. As for the penalized state, the stability condition was related to a partition of the transformed Hessian matrix being negative definite in the vicinity of the optimal solution. In both the relaxed and the penalized states, the authors proposed patch tests that could be performed under uniform stress states to test the stability of certain combinations of displacement and density interpolations \citep[p.~214 \& 220]{Jog1996}.

Two noteworthy remarks on the work of \citet{Diaz1995} and \citet{Jog1996}; \textbf{(i)} in addition to the factors mentioned in the first paragraph of this subsection, both studies agreed that an additional factor that influences the appearance of checkerboarding is the strain state of the patch of elements, and \textbf{(ii)} in contradiction to a commonly held misconception, higher order elements (namely eight-node and nine-node quads) don't completely resolve the checkerboarding problem for the element-wise constant implementation. Instead, there is still a possibility of checkerboarding with higher order elements, albeit it's much harder to induce and is mainly problem-dependent. \citet{Jog1996} even reported the appearance of checkerboarding with eight- and nine-node quads in the relaxed rank-2 laminates formulation.

Fortunately, any method that attempts to enforce mesh-independent designs would inherently get rid of the checkerboarding problem \citep[p.~71]{Sigmund1998}. However, since the non-penalized variable thickness sheet problem or the relaxed homogenized formulation using optimal microstructures do not suffer from the other mesh-dependence anomalies\footnote{\, Mainly because the non-convexity effects induced by domain discretization are remedied by the relaxation of the density design variables.} \citep[p.~266]{haber1994variable}, direct measures have to be implemented to resolve the checkerboarding problem.

It's worth noting that one additional technique of discouraging checkerboarding is the use of multiple finite elements per design element \citep{Zhou2001a}. This is somehow equivalent to using higher order shape functions in that both reduce the stiffness of the design element. The opposite of this technique; using multiple design elements per finite element, is termed multiresolution topology optimization. Such technique still results in internal checkerboarding similar to QR patterns within each finite element, which could be remedied using clever filtering techniques \citep{Gupta2018}.

\par
\vspace{1em}
Returning to our discussion on the variable thickness sheet problem, in order to avoid the use of \say{artificial} restriction methods that affect the feasible solution space (i.e., filtering, perimeter control, or gradient control methods), using higher order elements is more of a \say{natural} way to discourage checkerboarding, at least for academic if not for practical applications. Nonetheless, using higher-order elements in a partially-relaxed (i.e., penalized) state doesn't resolve mesh-dependency problems. In addition to the analytical convergence and uniqueness proof of the variable thickness sheet problem discussed in \citep{Petersson1999}, its convexity can - to a large extent - be proven numerically (cf. \citet[p.~298]{Cardoso2003}).

Even though the variable thickness sheet problem results in a nonpractical solution, it's still of some theoretical value. The popularity of the variable thickness sheet problem stems from two characteristics: \textbf{(i)} the fact that it is considered as the \textit{maximum bound on structural efficiency} for a certain domain discretization. In other words, material utilization is maximized as the solution space is expanded with the freedom of utilizing intermediate density elements. Once the solution space is restricted to enforce feasible discrete 0/1 solutions, the maximum achievable structural efficiency starts to decline, and \textbf{(ii)} it can be used to obtain \textit{the global relaxed optimum} which is the best candidate to be used as the initial guess towards discrete topology optimization.

In the following section, we discuss the different forms of penalization in use and how it affects the shape of the objective function and the global optimum.

\section{Penalization and its Effect on Convexity}
\label{sec_penal}
In a topology optimization context with a minimization objective, the essence of penalization is modifying the objective function in order to increase its value at certain \say{undesirable} design points in an effort to force the local minima towards the remaining unaffected \say{desirable} design points. For instance in a structural mechanics problem, the undesirable design points could be the ones of intermediate density, so penalization would be introduced to raise the objective function value at these points, whilst leaving the remaining desirable design points of discrete density unaffected. In topology optimization, penalization is utilized for one of two reasons: \textbf{(i)} Steering the solution towards a discrete design point to reverse the effect of relaxing the design variables as discussed in Section \ref{sec_var_thick_sh_prob}, or \textbf{(ii)} Enforcing a design constraint such as minimum/maximum length scales (e.g., manufacturability), and under this category also falls dealing with numerical instabilities such as checkerboarding, mesh dependency, and one-node hinges.

\subsection{Penalization of Intermediate Density Elements}
Penalizing intermediate density elements can be implemented in a number of different forms, they can be categorized as follows:

\textbf{\textit{Implicit --}} in which the penalization factor is not explicitly recognized or controlled. This type of penalization is achieved using homogenization with a sub-optimal microstructure (i.e., partial relaxation), that is a microstructure that is \textit{incapable} of obtaining the optimal energy bounds for intermediate density elements, hence penalization is inferred. A prominent example of this category is found in \citep{Bendsoe1988} where rectangular elements with rectangular holes were utilized (cf. \cite{Bendsoe1999} for a comprehensive review).
Historically, to the best of the authors' knowledge, this form of penalization was the precursor to modern density-based methods, as in the pioneering work by \citet{Bendsoe1988} which was followed by the first appearance of the SIMP method in the seminal work by \citet{Bendsoe1989}.
A noteworthy remark about this category is that the degree of penalization depends on the type of microstructure utilized in performing the homogenization.
Since the penalization factor cannot be explicitly controlled, it's often not possible to produce close-to-discrete designs using these approaches. Hence, such \textit{implicit} penalization should be combined with \textit{explicit global} penalization (see next) to steer the converged solution towards more discrete designs as was done in \citep{Allaire1993a}.
    
\textbf{\textit{Explicit --}} in which the penalization factor can be explicitly recognized and controlled, which is probably the more intuitive approach. This category includes two main approaches; \textbf{(i)} \textit{Local}; where the material property of interest (e.g., stiffness in a structural mechanics problem, permeability in a fluid flow problem, or conductivity in a heat transfer problem) is penalized explicitly within the local stiffness matrix of each element. This approach is more popular and is the basis of all density-based methods, and \textbf{(ii)} \textit{Global}; in which the objective function itself is appended with a penalty term that is only active for the undesirable intermediate density elements, and inactive otherwise.
This approach is often accompanied by complete relaxation of the design space either by using optimal microstructures such as rank-2 laminates (cf. \citet[p.~243]{Allaire1993}) or linear material interpolation such as the variable thickness sheet problem. Both the explicit local and global approaches have similar effects on the objective function. Albeit, a major difference is in the form of the governing equations being solved by the finite element solver.

Referring to our discussion on convex functions at the beginning of Section \ref{sec_var_thick_sh_prob}, penalization introduces non-convexity to the problem by transforming the objective function's epigraph from a convex set into a non-convex one. The effect of penalization on an originally convex objective function can be visualized in Fig. \ref{f1_b}. The value of the objective function at the desirable points (represented by the tick marks on the $x$-axis) remains unchanged, while it increases at the undesirable points (represented by the values in-between the tick marks). The last statement is significant in the sense that for a valid comparison between converged solutions of different penalization formulations/parameters, either the penalization factor has to be similar, or the converged solution has to be approximated to the nearest unaffected desirable point. For instance, in the case of penalizing intermediate density elements as in all density-based approaches, approximations towards the nearest discrete point could be achieved using a simple thresholding technique that still satisfies the volume fraction constraint (cf. \citet[p.~1051]{Sigmund2013} for an example code).

\subsection{Penalization for Enforcing Design Constraints}
In this subsection, we discuss in detail a few techniques that belong, or at least are relevant, to penalization for enforcing design constraints, and their effect on the problem's convexity. A summary of the Lagrangian multipliers method for enforcing constraints is in order first.

In most, if not all, the solvers currently available for topology optimization, design constraints are always enforced using the Lagrangian multipliers method. This is true for the Optimality Criteria (OC) method \citep{Bendsoe2004} where the only constraint is usually the simple volume fraction, and for the MMA \citep{Svanberg1987}, GCMMA \citep{Svanberg1995}, CONLIN \citep{Fleury1989}, IPOPT \citep{Wachter2006}, and SNOPT \citep{Gill2005} methods where generally any constraint can be implemented\footnote{\, Any reference to \say{solvers} in this article is intended to mean the mathematical solvers used to solve the optimization problem itself (e.g., OC, MMA, GCMMA, etc.) and not mathematical solvers of systems of equations (e.g., CG, GMRES, MINRES, etc.).}.
It's commonly known that the essence of the Lagrangian multipliers method is to transform the constrained optimization problem into an unconstrained problem where we look for stationary or critical points. This is unlike the common fallacy that Lagrangian multipliers turn the constrained optimization problem into an unconstrained optimization problem. In fact, the transformed problem is minimized with respect to the original design variables and maximized with respect to the Lagrangian multipliers, which makes the optimum of the original objective function a stationary point in the Lagrangian function space and not a minimum \citep{Kalman2009}. Another way of looking at the Lagrangian multipliers method is that it adds a penalty term to the objective function which becomes active when the constraints are violated and inactive otherwise. Depending on the characteristics of the penalty term, any mathematical constraint added to the optimization problem may in fact cause some type of non-convexification and should be investigated.

\subsubsection{Volume Fraction Constraint}
\label{vf_constr}
Compliance minimization problems must feature a volume fraction constraint. It's usually formulated as a linear summation of elemental densities weighted by their corresponding elemental volumes. Even though the volume fraction constraint is typically implemented as an inequality as can be seen in Eq. \ref{eqn_prob_form}, such form is only necessary for ease of mathematical implementation and to account for any minor variations due to the discrete nature of the desired final solution. Due to the physical nature of the compliance minimization problem and assuming a well-tuned mathematical algorithm, this constraint is almost always active in all iterations. Hence, without loss of generality, we can limit our discussion to the equality form. In addition, to only consider the potential non-convexification effects of this constraint, let's assume that the original problem has a convex set of feasible solutions; that is the relaxed form of the discretized problem (e.g., a variable thickness sheet problem). It's a known fact that the only equality constraint to produce a convex set is a linear one \citep[p.~289]{Cooper1970}. Such a linear constraint would transform the feasible set of solutions from the \(R^N\) space to a hyperplane in \(R^{N-1}\) space (which is always convex), where \(N\) is the number of finite elements. Hence, although the volume fraction constraint is enforced using Lagrangian multipliers, it actually doesn't introduce non-convexification to the problem. This is mainly due to its \say{linear equality} nature, which results in a still convex set of feasible solutions.

\subsubsection{Other Constraints Enforced through Appending the Objective Function}
In this subsection, we direct our attention to more sophisticated design constraints enforced through appending a penalty term to the objective function, either through Lagrangian multipliers or through an explicit penalty term as in \citep[p.~268]{haber1994variable}. \par
A careful look at the nature of the design constraints enforced through penalty functions or Lagrangian multipliers should provide some insight into the nature of the non-convexification imposed on the objective function as a result. By nature, in this context, we mean the locality vs. globality of the design constraint. With some abuse of terminology, locality vs. globality of the constraint could be loosely described as the extent to which the enforced design constraints affect the individual design variables.

To further elaborate, on one end of the spectrum, a perfect example of an extremely local constraint is the discreteness of the design variables\footnote{\, Technically, explicit global penalization of intermediate density elements can be thought of as a type of constraint imposed through using the Lagrangian multipliers method.}. Such a constraint dictates specific requirements on each design variable; that is each and every design variable has to take a certain value from an extremely limited set \(\{0,1\}\). Without taking such extremely local measures, the discreteness constraint cannot be satisfied. Hence, such local constraints would impose strong non-convexification effects on the objective function. In other words, the objective function cannot move from one discrete design point to another discrete point without coming across a penalized peak in between, hence a huge number of local minima.

On the other end of the spectrum, a perfect example of an extremely global constraint is the volume fraction constraint. Such a constraint doesn't impose any local measures on the design variables. Instead, its penalization effects described in Subsection \ref{vf_constr} act like a boundary surrounding an internal feasible solution space (i.e., the hyperplane in \(R^{N-1}\)) that is still a convex set. Hence, any combination of elemental densities equal to the prescribed volume fraction (i.e., a design point on the hyperplane) would satisfy the constraint, without any local restrictions on the design variables.

This discussion could also be extended to stress-based design constraints. Due to the local nature of the stress constraints, they introduce heavy non-convexification to the objective function. Global measures to overcome the locality problem typically reduce the non-convexification effects, albeit at the cost of not completely satisfying the constraints \citep{le2010}. A noteworthy remark is that the above \say{locality vs. globality} discussion is inherently related to the physics of the problem, and consequently would always apply without any dependence on the mathematical solver utilized.

\subsubsection{Mesh-Independency Filtering}
Since the first appearance of mesh-independency filters (i.e., sensitivity filtering by \citet{Sigmund1994, Sigmund1997}), they have become extremely popular in the research community for their effectiveness, computational efficiency, and ease of implementation. Generally, filtering methods are implemented as an intermediate step in the cycle of the analysis step (i.e., solving the equilibrium equations and calculating the original sensitivities) and the optimization step (i.e., updating the design variables). They can be divided into two categories: sensitivity filtering and density filtering\footnote{\, While some references treat Heaviside filtering as separate from density filtering, we choose to categorize it as a type of density filtering as it operates on the elemental densities.}. Sensitivity (density) filters work through modifying each element's sensitivity (density) based on the sensitivities (densities) of neighbouring elements within a predefined filter radius. Although these filters were originally invented to prevent common numerical instabilities such as checkerboarding and mesh dependency, some filters' capabilities have been extended to implementing minimum and maximum length scales in both solid and void regions in addition to some physics-specific constraints such as preventing one-node hinges in compliant structural mechanisms. A worthy remark is that although the original filter proposed by \citet{Sigmund1994, Sigmund1997} was of the sensitivity type, most filters currently in use are of the density type. One of the reasons behind this is that it's relatively easier to conceptually link the parameters used in defining the density filter to the density results, which is not easily the case with sensitivity filtering. For a comprehensive review, see \citet{Lazarov2016} and references therein. 

It seems prevalent that almost all the work published on black/white enforcing filters recommends using continuation methods with the filtering schemes (cf. \citet[p.~409]{Sigmund2007}, \citet[p.~249]{Guest2004}, and \citet[p.~128]{Guest2009a}). In other words, there is a consensus in the literature that sophisticated filtering schemes (i.e., filters that perform any task beyond the simple prevention of checkerboarding and mesh dependency), if implemented without continuation, would generally cause either unstable behavior at worst, or cause convergence to local minima at best. Yet, little attention is usually given to attempting to investigate the effects of the filtering schemes on the convexity of the problem. A probable cause for this general \textit{attitude} is the fact that the severe localized penalization imposed by density-based methods on intermediate density elements is generally much stronger and tends to overshadow any other form of non-convexification caused by additional constraints.

In the following, we attempt to investigate the effect of filtering methods on the convexity of the problem. At this point, it's worthwhile to revisit some generalizations of convex functions, namely \textit{quasiconvex functions}. Following the treatment in \citep[p.~134]{Bazaraa2006}, a function \(f:S \rightarrow R\), where \(S\) is a nonempty convex set in \(R^N\), is said to be \textit{quasiconvex} if for each \(\mathbf{x}_1\) and \(\mathbf{x}_2 \in S\), the following inequality is true:
\begin{figure*}
    \centering
    \subfloat[]{\label{f_quasic_a} \includegraphics[width=0.32\textwidth]{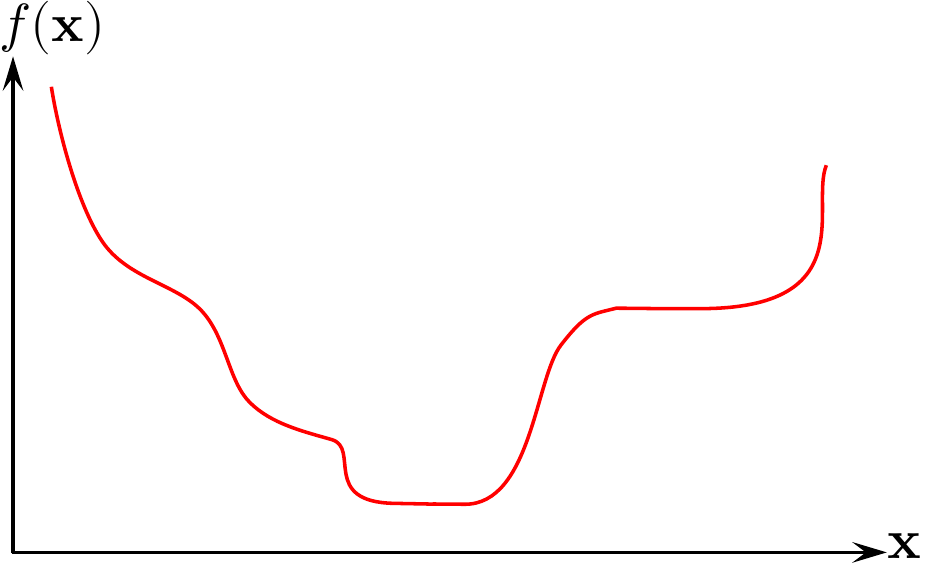}}
    \subfloat[]{\label{f_quasic_b} \includegraphics[width=0.32\textwidth]{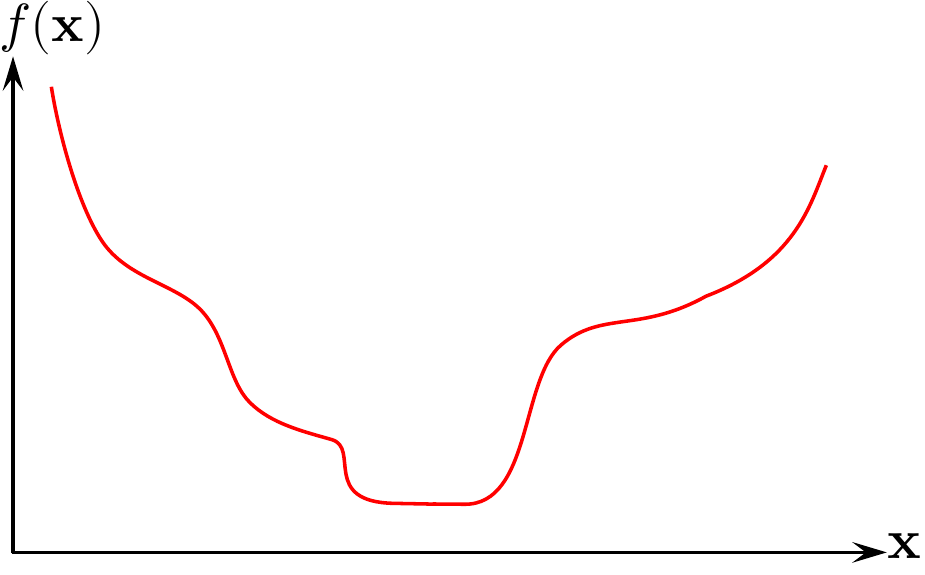}}
    \subfloat[]{\label{f_quasic_c} \includegraphics[width=0.32\textwidth]{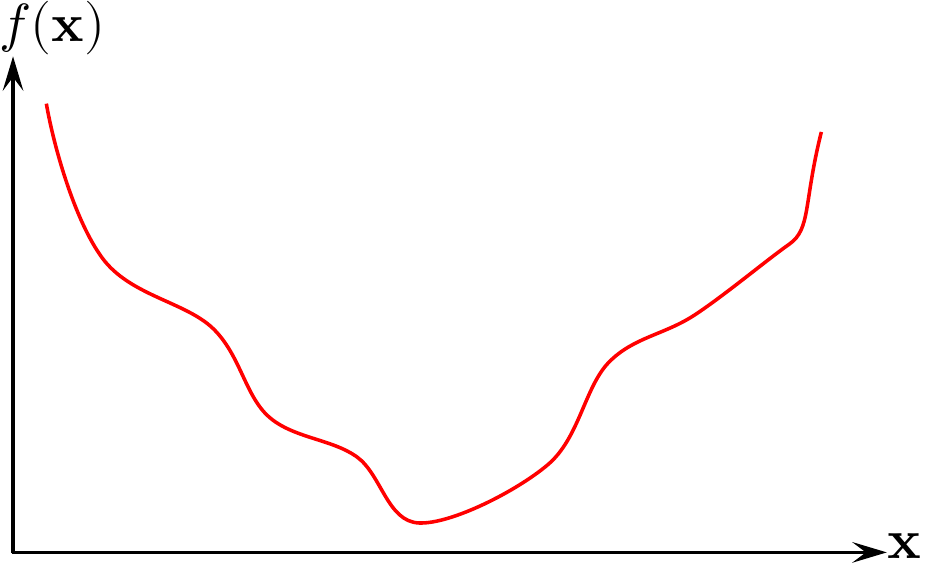}}
    \caption{Generalizations of convex functions: \textbf{a.} quasiconvex, \textbf{b.} strictly quasiconvex, and \textbf{c.} strongly quasiconvex. Notice the difference between the stationary points in the three types.}
    \label{f_quasi_convex}
\end{figure*}

\begin{equation}
\begin{aligned}
f[\lambda \mathbf{x}_1 + (1-\lambda) \mathbf{x}_2] \leq \text{max} \{f(\mathbf{x}_1), f(\mathbf{x}_2)\}, \\
\text{for each } \lambda \in (0,1).
\end{aligned}
\end{equation}

Although quasiconvex functions have a single minimum, their definition doesn't preclude the existence of multiple stationary points (see Fig. \ref{f_quasic_a}), which is detrimental to iterative solvers. Hence, a more useful concept is the strictly quasiconvex functions. A function \(f:S \rightarrow R\), where \(S\) is a nonempty convex set in \(R^N\), is said to be \textit{strictly quasiconvex} if for each \(\mathbf{x}_1\) and \(\mathbf{x}_2 \in S\) with \(f(\mathbf{x}_1) \neq f(\mathbf{x}_2)\), the following inequality is true:
\begin{equation}
\begin{aligned}
f[\lambda \mathbf{x}_1 + (1-\lambda) \mathbf{x}_2] < \text{max} \{f(\mathbf{x}_1), f(\mathbf{x}_2)\}, \\
\text{for each } \lambda \in (0,1).
\end{aligned}
\end{equation}

By enforcing \(f(\mathbf{x}_1) \neq f(\mathbf{x}_2)\) and the strictness of the inequality, we actively eliminate any stationary points except at the global minimum (see Fig. \ref{f_quasic_b}). If one is interested in a generalization of convex functions that supports uniqueness of solutions, \textit{strongly quasiconvex functions} might be of interest. A function \(f:S \rightarrow R\), where \(S\) is a nonempty convex set in \(R^N\), is said to be \textit{strongly quasiconvex} if for each \(\mathbf{x}_1\) and \(\mathbf{x}_2 \in S\), with \(\mathbf{x}_1 \neq \mathbf{x}_2\), the following inequality is true:
\begin{equation}
\begin{aligned}
f[\lambda \mathbf{x}_1 + (1-\lambda) \mathbf{x}_2] < \text{max} \{f(\mathbf{x}_1), f(\mathbf{x}_2)\}, \\
\text{for each } \lambda \in (0,1).
\end{aligned}
\end{equation}

\begin{figure}
    \centering
    \includegraphics[width=\columnwidth]{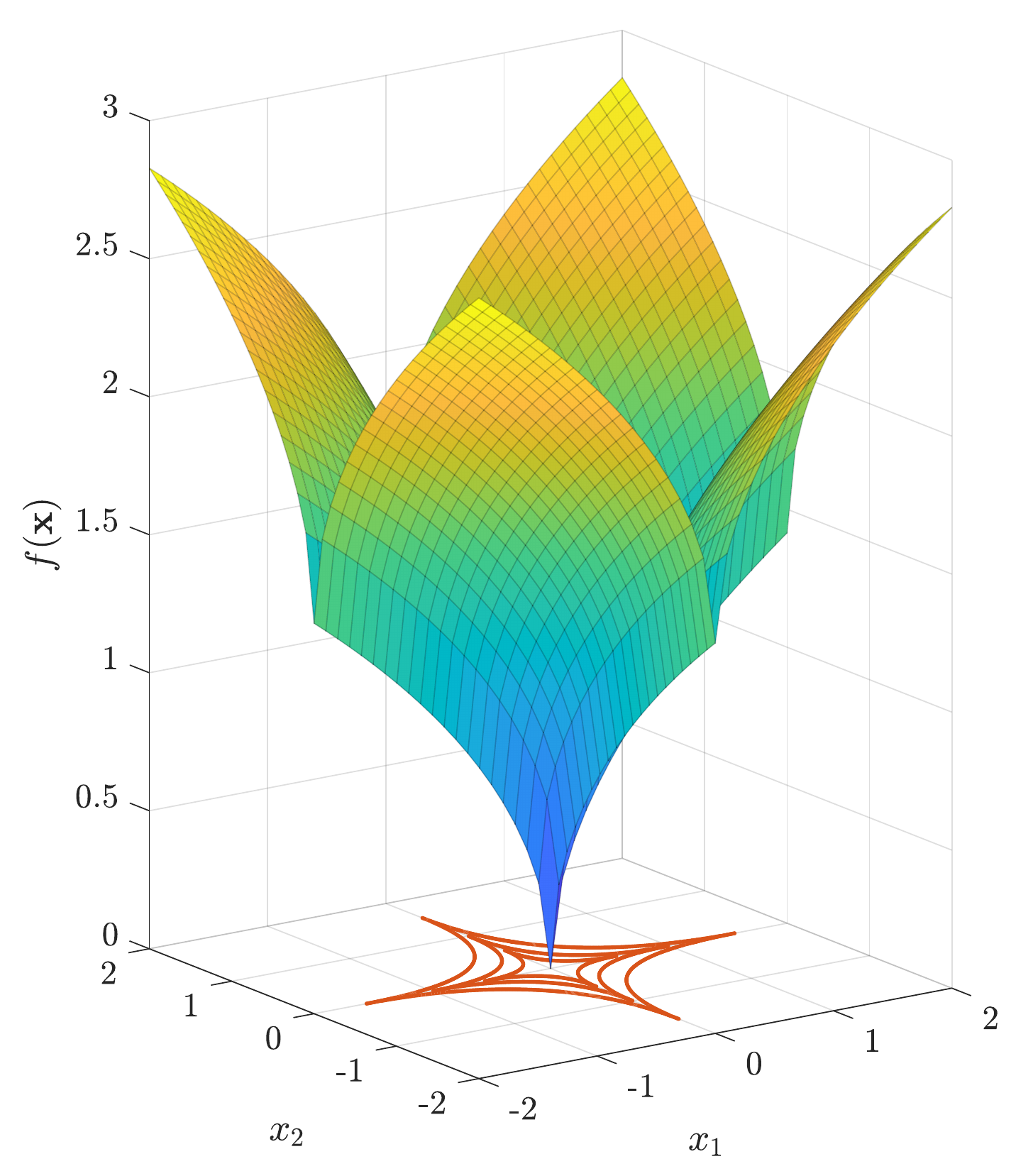}
    \caption{Function \(f(x_1, x_2)=\sqrt{|x_1|} + \sqrt{|x_2|}\) is an example of a multivariate unimodal function that doesn't belong to any category of quasiconvex functions since its lower level sets aren't convex. Yet, it has a unique global minimum that can be found using iterative solvers.}
    \label{f_unimodal}
\end{figure}

By enforcing \(\mathbf{x}_1 \neq \mathbf{x}_2\), strongly quasiconvex functions (see Fig. \ref{f_quasic_c}) assert uniqueness of the global optimum. However, a major disadvantage of the above mentioned functions is that they must have convex lower level sets. That is, for \(S\) a nonempty convex set in \(R^N\), the lower level sets defined as: 
\begin{equation}
\begin{aligned}
    L_\alpha = \{x \in S \, | \, f(x) \leq  \alpha\}.
\end{aligned}
\end{equation}

\noindent for \(\alpha \in R\), must be convex \citep[p.~95]{boyd2009}. This disadvantageous property is usually hard to prove for multivariate functions and doesn't add any significant value for an iterative solver. Hence, it would be useful to explore a more generalized concept, that is \textit{unimodal functions}. \citet[p.~156]{Bazaraa2006} define \textit{univariate} unimodal functions as follows; a function \(f : S \rightarrow R\), where \(S \text{ is some interval on } R\), is unimodal on \(S\) if there exists an \(x^* \in S\) at which \(f\) attains a minimum and \(f\) is nondecreasing on the interval \(\{x \in S: x^* \leq x\}\) and nonincreasing on the interval \(\{x \in S: x \leq x^*\}\). As for \textit{multivariate} unimodal functions, a number of conflicting definitions exist in literature, specially in the field of Probability and Statistics \citep{Dharmadhikari1976, Perlman1988, Dharmadhikari1988}. Some authors associate the definition of unimodality with that of quasiconvexity as in having convex lower level sets, which is true for univariate functions but not for multivariate ones. In this work, for the definition of multivariate unimodal functions, the natural extension of univariate unimodal functions is utilized \citep[p.~38]{Dharmadhikari1988}. That is; a multivariate function \(f(\mathbf{x})\) is called unimodal if it's nondecreasing along rays emanating from its global minimum in all directions (see Fig. \ref{f_unimodal}). To preclude the existence of multiple stationary points, we could utilize the notion of strict monotonicity so that \(f(\mathbf{x})\) is strictly increasing in every direction emanating from its global minimum.

Returning to our discussion on the effect of filtering on convexity, it's true that convexity might be distorted. However, unimodal functions are an example of non-convex functions that still possess a global minimum attainable through iterative solvers. Hence, it's safe to assume that as long as the modified function is strictly unimodal, no new local minima would arise because of this non-convexity. It's worth noting that in penalized topology optimization, the unfiltered problem is already non-convex with many local minima. Hence, the above discussion would be only applicable to smaller \textit{locally} convex parts of the large non-convex function (i.e., a localized valley that has a single minimum). In the following, we reflect this discussion onto sensitivity filtering.

\textit{i -- Sensitivity Filtering}
What the sensitivity filtering does is that it modifies the descent direction used in the optimizer, albeit without affecting the value of the objective function itself \citep[p.~474]{Sigmund2012a}. One could say that sensitivity filtering \say{tricks} the optimizer by altering the gradient of the objective function. Hence, the optimizer would follow a different route towards the minimum and obviously would settle at a different minimum from the unfiltered original minimum, mainly because the sensitivity filtering flattened the objective function at this point as seen from the gradient's perspective. Since the gradient has been altered, it's safe to say that sensitivity filtering \say{might} cause some degree of non-convexification, albeit it cannot strictly be called penalization since the value of the objective function itself hasn't been altered.

Any filtering method starts by determining the neighbourhood of each element; that is the set \(N_e\) consisting of the elements with centers spatially located within a given filter radius \(r\) of the center of element \(e\) as follows:
\begin{equation}
    \begin{aligned}
        N_e = \{i \; | \; || \mathbf{x}_i - \mathbf{x}_e || \leq r \}. \\
    \end{aligned}
\end{equation}
\noindent where \(\mathbf{x}_i\) denotes the spatial location of the center of element \(i\).

A modified version of the sensitivity filter that accounts for non-regular meshes with varying elemental volumes is as follows \citep[p.~408]{Sigmund2007}:
\begin{subequations}
\begin{align}
& \widetilde{\frac{\partial f}{\partial \rho_e}} = %
\frac{\mathlarger{\sum}\limits_{i \in N_e} w_e(\mathbf{x}_i) \rho_i \mathlarger{\frac{\partial f}{\rho_i}}/v_i} %
{\rho_e/v_e \mathlarger{\sum}\limits_{i \in N_e} w_e(\mathbf{x}_i)} , \label{eqn_sens_filt} \\[4pt] 
& w_e(\mathbf{x}_i) = r - || \mathbf{x}_i - \mathbf{x}_e || .
\label{eqn_linear_weight}
\end{align}
\end{subequations}

\noindent where \(f\) is the objective function, \(\rho_e\) is the density of element \(e\), \(w_e(\mathbf{x}_i)\) is a linearly decaying (cone-shaped) weighting function\footnote{\, The upcoming proof still applies to other weighting functions such as the Gaussian (bell-shaped) distribution and the constant weighting functions since they are always positive.} of element \(i\) within the neighbourhood of element \(e\), \(v_i\) is the volume of element \(i\), and \(\partial f/\partial \rho_i \; \text{and} \; \widetilde{\partial f/\partial \rho_i}\) denote the original and modified sensitivity respectively of the objective function with respect to element \(i\).

A noteworthy remark is on how the descent direction is calculated inside the typical mathematical solvers used in topology optimization (e.g., OC, GCMMA, MMA, etc.). The general concept is that the solvers try to satisfy the Karush-Kuhn-Tucker conditions by utilizing a combination of the gradients of the objective function and the constraints. Subsequently, the descent direction and step are calculated based on the given move limits and the bounds of the design variables. Consequently, it is safe to say that as long as the inputted gradients denote ascent directions, the outputted updated design variables would always denote a reduction in the objective function (except at a minimum of course). Hence, our investigative approach for sensitivity filters would be focused on the input to the mathematical solvers; that is whether the filtered sensitivities still constitutes an ascent direction or not. In other words, whether the modified function is \textit{locally} strictly unimodal or not. In mathematical terms, it means satisfying the following condition:
\begin{equation}
\begin{aligned}
\boldsymbol{\frac{\partial f}{\partial \rho}} \cdot \boldsymbol{\widetilde{\frac{\partial f}{\partial \rho}}} > 0.
\label{equation9}
\end{aligned}
\end{equation}

\noindent In order to ensure that such a dot product is always positive, we need to ensure that each multiplication term is actually positive on its own as follows:
\begin{equation}
\begin{aligned}
\frac{\mathlarger{\sum\limits_{i \in N_e}} w_e(\mathbf{x}_i) \rho_i \mathlarger{\frac{\partial f}{\rho_i}}/v_i}{\rho_e/v_e \mathlarger{\sum\limits_{i \in N_e}} w_e(\mathbf{x}_i)} \cdot \frac{\partial f}{\partial \rho_e} > 0.
\label{eqn_filt_sens_mulp}
\end{aligned}
\end{equation}

Since \(w_e(\mathbf{x}_i), \, \rho_i, \, \text{and} \, v_i\) are always positive, any multiplicative combination of these quantities would always be positive. Given the fact that sensitivities are always negative in compliance minimization problems, this means that the modified sensitivities are always negative. Hence, each individual term in Eq. \ref{eqn_filt_sens_mulp} is in fact positive as it's a multiplication of two negative quantities. This concludes the proof that the filtered sensitivities in compliance minimization problems always constitute an ascent direction (i.e., the modified function is \textit{locally} strictly unimodal). A worthy remark is that an essential component of our proof is that in compliance minimization problems, the original sensitivities are always negative. Hence, in other problems where the sensitivities might be either positive or negative, there exists the possibility that sensitivity filtering might cause non-convexification and disturb the trajectory of the mathematical solver.

\textit{ii -- Basic Density Filtering}

This category includes types of density filtering that enforce a grey transition region along the boundaries. Density filtering works by mapping each design point to another design point based on the details of the density filtering scheme. This mapping results in two main effects; \textbf{(i)} a jump from the original design point to the filtered one, and \textbf{(ii)} based on this jump, the sensitivity has to be modified. A typical density filtering takes the following form \citep{Bruns2001a, Bourdin2001}:
\begin{equation}
\begin{aligned}
\widetilde{\rho_e} = \frac{\mathlarger{\sum\limits_{i \in N_e}} w_e(\mathbf{x}_i) v_i \rho_i}{\mathlarger{\sum\limits_{i \in N_e}} w_e(\mathbf{x}_i) v_i}.
\label{eq_dens_filt}
\end{aligned}
\end{equation}

\noindent and the sensitivities are modified accordingly as follows:
\begin{subequations}
\begin{align}
\label{eq_dens_filt_sens_a}
& \frac{\partial f}{\partial \rho_e} = \sum\limits_{i \in N_e} \frac{\partial f}{\partial \widetilde{\rho_i}} \frac{\partial \widetilde{\rho_i}}{\partial \rho_e}, \\[4pt]
\label{eq_dens_filt_sens_b}
& \frac{\partial \widetilde{\rho_i}}{\partial \rho_e} = \frac{w_e(\mathbf{x}_e) v_e}{\mathlarger{\sum\limits_{j \in N_i}}w_e(\mathbf{x}_j) v_j}. 
\end{align}
\end{subequations}

Let's focus on the resulting sensitivities first, it's clear that \(\partial \widetilde{\rho_i} / \partial \rho_e\) is always positive since it's a multiplicative combination of \(w_e(\mathbf{x}_e)\) and \(v_e\) which are always positive. Hence, according to Eq. \ref{equation9}, the modified sensitivities (Eq. \ref{eq_dens_filt_sens_a}) didn't change signs and still in fact constitute an ascent direction. This proof, unlike the sensitivity filtering proof, applies to any problem, not just to compliance minimization.

As for the modified densities, it's worthy to investigate whether the filter constitutes an affine transformation or not. A careful look at the modified densities (Eq. \ref{eq_dens_filt}) reveals that it can be put in the form:
\begin{equation}
\begin{aligned}
\{ \bm{\widetilde{\uprho}} \} = \lbrack \mathbf{A} \rbrack \{\bm{\uprho}\}.
\end{aligned}
\end{equation}

\noindent where \(\mathbf{A}\) is an \(N \times N\) square matrix since \(\bm{\widetilde{\uprho}}\) and \(\bm{\uprho}\) have the same dimension \(N\). This relation constitutes a linear mapping (or an affine mapping to be more general \citep[p.~79]{boyd2009}) from the original to the modified densities. A rather nice property of affine transformations is that the resultant set is convex if the original is, i.e., the modified densities set is convex if the original densities set is convex, which it is. Relating to our discussion on convex functions in Section \ref{sec_var_thick_sh_prob}, this concludes that the modified domain is indeed convex. 
As for the convexity of the codomain, since the objective function is merely evaluated at a new design point, its codomain undergoes no changes and retains its original convexity, if any existed.

Although we have proven the basic density filter maintains convexity, an intriguing question remains; does the affine mapping exclude any design points from the modified domain? In order to answer this question, we need to investigate whether the matrix \(\mathbf{A}\) is invertible or not. This is mainly because a non-invertible affine mapping collapses the space along some directions, i.e., the modified set is smaller than the original one. The matrix \(\mathbf{A}\) is singular if its determinant is zero, which can happen in one of two cases (or both): \textbf{(i)} a whole row/column is zero, or \textbf{(ii)} two rows/columns are linearly dependent. The first case cannot happen since the worst case scenario - i.e., having a filtering radius \(r\) that doesn't span any adjacent elements - will result in an identity matrix. Any increase in the filtering radius \(r\) would fill in more zero elements and hence it's impossible for \(\mathbf{A}\) to have a whole row/column of zeros. As for the second case, the only scenario that would result in two (or more) rows/columns being linearly dependent is if two elements have the same neighbourhood of elements \(N_e\) and a constant weighting function \(w_e(\mathbf{x}_i)\). The only case we can see this happening is at the free edge of a thin structure comprised of two rows of finite elements, the two elements at the edge would have the same neighbourhood of elements. Of course this is an extremely unpractical\footnote{\, Even with including the elements outside boundary as void elements, the two elements could still have the same neighbourhood of elements (cf. \citep[p.~407]{Sigmund2007} for more details on the treatment of mesh boundaries).} case and could be safely ignored. Hence, for all practical purposes, the modified domain of a basic density filter doesn't exclude any design points. We conjecture that a careful choice of a weighting function combined with a specific design point, though unrealistic, could result in a checkerboard pattern as a result of filtering. It goes without saying that this is purely a mathematical observation and doesn't happen in any practical situation.

\textit{iii -- Filtering that Enforces 0/1 Designs}

As mentioned before, sensitivity as well as basic density filters enforce a grey transition region along the boundaries. To overcome this issue, a class of filters that enforce 0/1 discrete designs were developed. The first appearance of such filters was in the pioneering work by \citet{Guest2004}, in which the authors used a nonlinear projection (a regularized Heaviside step function) to ensure a discrete 0/1 boundary and enforce a minimum length scale on the solid phase. A slight modification of this filter could enforce the minimum length scale on the void instead of the solid phase, cf. \citet[p.~125]{Guest2009a}.

It's well known that intermediate density penalization is mainly achieved by the penalization of the objective function itself, and is usually enforced gradually through a continuation method. To get rid of the grey transition regions, these filtering schemes enforce the jump to be towards a discrete rather than an intermediate density design point, which can be considered a form of penalization. Intermediate density penalization, if enforced by the filtering scheme rather than the objective function, would cause unstable behavior (cf. \citet[p.~249]{Guest2004} and \citet[p.~409]{Sigmund2007}). Mainly because the filtering scheme is not designed to seek the minimum discrete point of the problem and consequently its modified sensitivities wouldn't lead to that direction. Hence, it's essential that continuation is utilized in the filtering scheme so as to follow a similar (or lower) degree of penalization as that of the objective function. In what follows, a mathematical justification is presented.

\citet[p.~248]{Guest2004}'s Heaviside filter takes the following form:
\begin{equation}
\begin{aligned}
\Bar{\rho}_e = 1 - e^{-\beta \widetilde{\rho}_e} + \widetilde{\rho}_e \, e^{-\beta}.
\label{eqn_heav}
\end{aligned}
\end{equation}

\noindent where the parameter \(\beta\) controls the curvature of the regularization (\(\beta = 0\) recovers the basic density filter and \(\beta \to \infty\) recovers the Heaviside step function), and \(\widetilde{\rho}\) is calculated as in Eq. \ref{eq_dens_filt} using a linear weighting function as in Eq. \ref{eqn_linear_weight} (other weighting functions could also be used). The affine mapping argument used in proving the convexity of basic density filtering is not valid here since its converse is not true (i.e., not being an affine mapping doesn't necessarily mean producing a non-convex set). Hence, a different argument is needed, namely the generalization of which the affine mapping argument is a special case.

\citet[p.~114]{soltan2015lectures} stated that \say{a mapping \(f : R^N \to R^M\) is called \textit{convexity-preserving} if the \textit{f}-images of all convex sets in \(R^N\) are convex sets in \(R^M\)}. In our case, \(M = N\). Hence, we can focus our attention on a single modified density \(\Bar{\rho}\), whether its output is a convex set or not. It's clear from the analytical form of the Heaviside filter (recovered as \(\beta \to \infty\)) that the result contains only two elements \(\{0,1\}\). In other words, any combination of neighbourhood elements' weights and densities would always result in a filtered density of 0 or 1. This means that certain design points would be excluded from the modified set, and in our case the excluded points are the intermediate density ones. Hence, it's clear that the result is not a convex set, and the filter is not a convexity-preserving mapping. It's worth noting the resemblance between the effect of these filters and the discussion on domain discretization in Section \ref{sec_var_thick_sh_prob}.


\textit{iv -- Other Filters}

This subsection concerns the filters that don't enforce neither grey nor discrete designs. To overcome the downsides of the Heaviside filter, \citet{Sigmund2007} introduced the morphology-based operators \textit{dilate} and \textit{erode}, which enforce a minimum length scale on the solid and void phases respectively. These filters work as maximum (\textit{dilate}) or  minimum (\textit{erode}) operators, meaning they will select the maximum (or minimum) value of the neighbourhood elements\footnote{\, Assuming constant weighting as in \citet{Sigmund2007}'s original formulation. Using a different weighting function will slightly alter this definition.}. Later, \citet{Svanberg2013} provided alternative definitions of the \textit{dilate/erode} operators based on the second and third Pythagorean means (i.e., geometric and harmonic means) instead of the usual first Pythagorean mean (i.e., arithmetic mean).

It's worth noting that these filters are nonlinear (not affine mappings), but they can produce a convex set for each \(\Bar{\rho}\). With various combinations of neighbourhood elements, they can produce \(\Bar{\rho} \in [0,1]\). In fact, \citet[p.~865]{Svanberg2013} proved that the three dilate operators (i.e., morphology, geometric, and harmonic) are convex density filters (i.e., each \(\Bar{\rho}\) is a convex function of \(\bm{\uprho}\)) through analytically proving that the Hessian is always positive semidefinite.

This concludes our discussion on the convexity of filters, we next move to discussing the various continuation parameters in use.

\subsection{Continuation Methods}
\label{sec_contn}
In a topology optimization context, continuation refers to the process of gradually increasing the degree of non-convexity of the objective function. This non-convexity could be caused by a number of sources, albeit they all take the shape of penalization; either implicitly or explicitly.

To the best of the authors' knowledge, the first work to use continuation methods in a modern topology optimization context was that by \citet[p.~211]{Allaire1993a} and the follow-up work by \citet{Allaire1993}. Their work used the early method of homogenization using ranked laminates (i.e., different laminates on different microscales, cf. \citet[p.~165]{kohn1988recent} for details) to introduce controlled porosity into the structure. Since ranked laminates are considered an optimal microstructure (i.e., they provide efficient utilization of materials and voids) such formulation produces considerable areas of perforated structures. To approach the results of other less-optimal structures for a valid visual comparison, \citet{Allaire1993} introduced penalization of composite areas. To overcome the non-convexification effects caused by the penalization, the authors first computed the optimal design of the relaxed formulation, then they introduced a constant penalization factor starting from this relaxed optimal design. It's worthy to mention that, in \citet{Allaire1993a} and \citet{Allaire1993} due to the nature of homogenization, the authors had to use a global explicit penalization approach, where a penalty term is directly affixed to the objective function in such a way that it is active for intermediate density elements, and inactive otherwise. In addition, they recommended optimizing the problem without penalization (i.e., the relaxed version) first, then using a constant penalization factor from there on.

Another early mention of continuation methods is in the work of \citet{Haber1996}. The authors added an external penalization term to the objective function and gradually increased the penalization parameter in successive iterations. This implementation was performed in the process of producing sample results to be compared to the results of the main thesis of their work which was the parameter constraint method, that didn't include continuation. In another work, \citet[p.~1429]{Petersson1998} recommended starting with unity penalization till convergence, then increasing the penalization factor by \(0.5\) with convergence after each step. More recently, a detailed study on the effect of the continuation procedure on the objective function was performed by \citet{Rojas-Labanda2015a} in addition to a detailed literature survey. The current consensus is that continuation methods have to be used, with little emphasis on the step or the convergence criteria used in updating the penalization parameter.

\section{Challenges and Recent Trends}
\label{sec_challenges}
In this section, we discuss a few key concepts that are not directly related to topology optimization but are still relevant to non-convexity and its associated local minima problem. As detailed in the previous sections, in any optimization problem, the main detrimental outcome of non-convexity is the existence of local minima. In essence, the problem of local minima arises mainly with computational methods that follow a strict gradient-based trajectory starting from a specific initial guess. Probably the only solution to this problem in the case of gradient-based methods is starting from the minimum of an approximate relaxed, convex problem and using continuation methods to approach the global minimum as much as possible. Non gradient-based methods could actually be of use in this case since they don't follow a strict trajectory determined by the gradients and hence are less prone to \say{getting stuck} in a local minimum \citep{Huang2010}. For instance, \citet{Fujii2018} demonstrated the robustness of their topology optimization method based on a covariance matrix adaptation evolution strategy against the choice of the initial guess. One of the main arguments against using evolutionary methods is that it tends to take more iterations to arrive at a solution, however, recent works on approximate reanalysis of structures could be of help \citep{Senne2019}.

Global optimization offers an alternative to typical local optimization where a global search/exploration of the design space is performed some way or another \citep{Chunna2020}. Global optimization can also be combined with parallel computing to better manage the computational expenses \citep{Xing2020}. Certain probabilistic methods that have the capability of utilizing multiple starting points such as simulated annealing is also a potential player to overcome the local minima problem \citep{Zeng2019}. Finally, neural networks show potential in providing robust results even with relatively bad initial guesses \citep{swirszcz2017local}.

\section{Conclusions}
\label{sec_conc}
Convexity is a core concept in optimization given the abundance of convex programming algorithms and the fact that global optimality can be proven easily in convex problems. However, the current formulation of topology optimization problems inherently introduces non-convexification into the problem. This non-convexification is introduced first by the simple domain discretization required to solve the problem numerically, and later by the intermediate density penalization required to enforce discrete solutions. Various design constraints can be sources of additional non-convexity. In this article, we presented a comprehensive treatment of the various sources of non-convexification introduced into a topology optimization problem. The discretization of the design domain represents a restriction method with a huge pool of local minima. This non-convexity is then remedied through a relaxation of the design variables. In addition, we provided a mathematical treatment on effects of some of the common filtering methods on the convexity of the problem. In conclusion, the non-convexification effects of the penalization of intermediate density elements simply overshadows any other type of non-convexification introduced into the problem, mainly due to its severity and locality. Continuation methods are strongly recommended to overcome the problem of local minima, albeit its step and convergence criteria are left to the user depending on the type of application. 

\bibliographystyle{elsarticle-harv}
\pdfbookmark[0]{References}{References}
\bibliography{library}

\end{document}